\documentclass[ECP]{ejpecp}

\usepackage{amsmath,amssymb,amsfonts,amsthm,graphicx}
\usepackage[numbers]{natbib} % For \citep and \citet
\usepackage{mathtools} % For \vcentcolon to define : = and = :
\usepackage{xcolor} % To color the \fred and \CG comments
\usepackage{dsfont} % for \ind (indicator functions)

\SHORTTITLE{Stein's method for the symmetric matrix normal distribution}

\TITLE{Stein's method for the symmetric matrix normal distribution with an application to the approximation of the Wishart law \support{R.G.\ is supported
 by the Engineering and Physical Sciences Research Council (EP/Y008650/1) and F.O.\ is supported by the Natural Sciences and Engineering Research Council of Canada (RGPIN-2026-04471; DGECR-2026-00449).}}

\AUTHORS{%
 Robert~E.~Gaunt\footnote{The University of Manchester, United Kingdom. \EMAIL{robert.gaunt@manchester.ac.uk}} \and
 Fr\'ed\'eric~Ouimet\footnote{Universit\'e du Qu\'ebec \`a Trois-Rivi\`eres, Canada. \EMAIL{frederic.ouimet2@uqtr.ca}}}

\KEYWORDS{Stein's method; symmetric matrix normal; symmetric Ornstein--Uhlenbeck process; Wasserstein distance} % Separate items with ;

\AMSSUBJ{60F05} % Edit. Separate items with ;
\AMSSUBJSECONDARY{60H10; 60J60; 62E20; 62H10; 62H12} % Optional, separate items with ;

\SUBMITTED{June xx, 2026} % Edit.
\ACCEPTED{xx, 2026} % Edit.

\VOLUME{0}
\YEAR{2026}
\PAPERNUM{0}
\DOI{10.1214/YY-TN}

\ABSTRACT{In this paper, we extend Stein's method to the symmetric matrix normal distribution. In particular, we obtain a Stein characterization of the symmetric matrix normal distribution involving the extended generator of the symmetric matrix Ornstein--Uhlenbeck process, present a semigroup representation of the solution of the corresponding Stein equation, and establish regularity estimates for the solution. This framework of Stein's method for symmetric matrix normal approximation complements the recent theory of Stein's method for matrix normal approximation, and we make an explicit connection between these frameworks. We apply this theory to derive a Wasserstein distance bound for the symmetric matrix normal approximation of the Wishart distribution.}

\newcommand{\N}{\mathbb{N}}
\newcommand{\R}{\mathbb{R}}

\newcommand{\EE}{\mathsf{E}} % Russian style, do not change
\newcommand{\Var}{\mathsf{Var}} % Russian style, do not change
\newcommand{\Cov}{\mathsf{Cov}} % Russian style, do not change
\newcommand{\bb}[1]{\boldsymbol{#1}}
\newcommand{\rd}{\mathrm{d}}
\newcommand{\ind}{\mathds{1}}
\newcommand{\tr}{\mathrm{tr}}
\newcommand{\etr}{\mathrm{etr}}

\newcommand{\vecp}{\mathrm{vecp}}
\newcommand{\vecs}{\mathrm{vecs}}
\newcommand{\OO}{\mathcal{O}}

\newcommand{\leqdef}{\vcentcolon=}

\allowdisplaybreaks

\begin{document}

\section{Introduction}\label{sec:intro}

Stein's method is a powerful technique for bounding the distance between two probability distributions with respect to a probability metric. The method was introduced by \citet{s72} in the context of normal approximation, for which the theory is now very well developed; see the monographs \citep{MR2732624,np12} for introductions and examples of some of the numerous applications of Stein's method throughout the mathematical sciences. Stein's method was extended to the multivariate normal distribution by \citet{MR1035659} and \citet{g91}, among others. This theory is now rather sophisticated, as showcased by the recent applications to the quantification of the multivariate central limit theorem in the convex and Wasserstein distances \citep{bonis, MR3980309,MR4003566}.

Stein's method for multivariate normal approximation has found recent applications to distributional approximations arising in random matrix theory. In these applications, a vectorization argument is applied to transfer the analysis to a multivariate normal approximation \citep{m22,nz22,tudor}. However, vectorizing a random matrix by stacking its columns obscures the inherent algebraic structure, left-right symmetries, and intrinsic dependencies found in matrix models, meaning that in many settings one cannot transfer a matrix-variate distributional approximation to a multivariate approximation through column-stacking; see Section~4.2 of \citet{GauntOuimetRichards2026} for an example. This serves as motivation for a systematic development of Stein's method for matrix-variate distributions, which has recently been initiated by \citet{GauntOuimetRichards2026} and \citet{BaillyGauntOuimetRichardsvonSachs2026} for the matrix normal and Wishart distributions, respectively.

In this paper, we extend Stein's method to the symmetric matrix normal (SN) distribution.
The SN distribution is a symmetric-matrix analogue of the multivariate normal distribution and contains the Gaussian orthogonal ensemble (GOE) as a special case. Basic properties of the SN distribution can be found in the standard reference \citep[Section~2.5]{GuptaNagar2000}, and some application areas are collected in \citet{MR4358612}. In this paper, by making a connection to the framework of Stein's method for matrix normal approximation developed by \citet{GauntOuimetRichards2026}, we are able to efficiently establish the basic ingredients of Stein's method for SN approximation. We also apply this theory to derive a Wasserstein distance bound for the SN approximation of the Wishart distribution. Our bound complements a recent bound on the rate of convergence in total variation distance given by \citet[Theorem~6]{MR4358612} and cannot be derived via naive vectorization, thereby demonstrating the utility of the SN Stein framework.

The rest of the paper is organized as follows. Section~\ref{sec:definitions} collects notation and basic definitions. Section~\ref{sec:symmetric.matrix.normal.Stein.theory} develops Stein's method for the SN distribution. Section~\ref{sec:Wishart.approximation.shifted.symmetric.normal} states the Wasserstein distance bound for the SN approximation to the Wishart law. Section~\ref{sec:proofs} contains the proofs of all our results.

\section{Definitions and notation}\label{sec:definitions}

Throughout, $[d] \leqdef \{1, \ldots, d\}$ for $d\in \N \equiv \{1, 2, \ldots\}$. Let $\mathcal{S}^d$, $\mathcal{S}_{+}^d$, and $\mathcal{S}_{++}^d$ denote the sets of real symmetric, nonnegative definite, and positive definite $d\times d$ matrices, respectively. Unless otherwise stated, these spaces are equipped with the Frobenius inner product $\langle A, B \rangle_F = \tr(A^{\top} B)$ and the induced norm $\|A\|_F = \sqrt{\langle A,A\rangle_F}$. For any square matrix $A$, let $\tr(A)$ be its trace, $\etr(A) \leqdef \exp\{\tr(A)\}$, and $|A|$ its determinant. For $S\in \mathcal{S}_{+}^d$, the matrix $S^{1/2}$ denotes the symmetric square root and $\|S\|_2$ the spectral norm. If $B$ is an open subset of a finite-dimensional real Euclidean space, in particular if $B\subseteq \mathcal{S}^d$ is open in the relative topology, and $m\in \N_0 \equiv \{0,1,\ldots\}$, let $C^{m}(B)$ be the class of real-valued functions $f:B\to\R$ that are $m$ times continuously differentiable on $B$ (all partial derivatives up to total order $m$ exist and are continuous), and let $C_b^{m}(B)$ be the subclass for which all partial derivatives up to total order $m$, including the order-$0$ derivative $f$ itself, are bounded on $B$. For a map $F:B\to E$, where $E$ is a finite-dimensional real Euclidean space, $D^kF(x)$ denotes the $k$-th Fr\'echet derivative of $F$ at $x$, viewed as a $k$-linear map on the ambient real vector space of $B$ with values in $E$; its evaluation in directions $u_1,\ldots,u_k$ is written $D^kF(x)[u_1,\ldots,u_k]$ and, when $B\subseteq\mathcal{S}^d$, the directions $u_i$ lie in $\mathcal{S}^d$. The symbols $\bb{0}_d$, $0_{d\times d}$, and $I_d$ denote the $d$-dimensional zero vector, the $d\times d$ zero matrix, and the $d\times d$ identity, respectively.

Let $q_d \leqdef d(d + 1)/2$, and fix once and for all an ordering $\{(i_a,j_a):a\in [q_d]\}$ of the index set $\{(i,j):1\leq i\leq j\leq d\}$. Define the orthonormal half-vectorization
\[
\vecs(S) \leqdef (r_a S_{i_a j_a})_{a\in [q_d]}\in \R^{q_d}, \qquad
r_a \leqdef
\begin{cases}
1, &\mbox{if } i_a = j_a, \\
\sqrt{2}, &\mbox{if } i_a < j_a.
\end{cases}
\]
Then $\vecs:\mathcal{S}^d\to\R^{q_d}$ is an isometry, so $\|\vecs(S)\|_2 = \|S\|_F$. Let $\mathrm{mat}_{\mathrm{s}}:\R^{q_d}\to\mathcal{S}^d$ denote the inverse map, and define
\[
E_a \leqdef \mathrm{mat}_{\mathrm{s}}(\bb{e}_a), \qquad a\in [q_d],
\]
where $\bb{e}_1,\ldots,\bb{e}_{q_d}$ are the standard basis vectors of $\R^{q_d}$. Then $E_1,\ldots,E_{q_d}$ form an orthonormal basis of $\mathcal{S}^d$ for the Frobenius inner product. If $f:\mathcal{S}^d\to\R$, define
\[
\widetilde{f}(\bb{x}) \leqdef f(\mathrm{mat}_{\mathrm{s}}(\bb{x})), \qquad \bb{x}\in \R^{q_d}.
\]
For $\bb{a} = (a_1,\ldots,a_k)\in [q_d]^k$, define the coordinate directional derivative
\[
D_{\bb{a}}^{\mathrm{s}}f(S) \leqdef D^k f(S)[E_{a_1},\ldots,E_{a_k}], \qquad S\in \mathcal{S}^d,
\]
whenever the derivative exists. We use the convention $D_{\varnothing}^{\mathrm{s}}f = f$. Equivalently, $D_{\bb{a}}^{\mathrm{s}}f(S)$ is the corresponding coordinate partial derivative of $\smash{\widetilde{f}}$ at $\vecs(S)$. We use $\nabla f$ for the symmetric gradient, characterized by
\[
Df(S)[U] = \tr\{U\nabla f(S)\}, \qquad S,U\in \mathcal{S}^d.
\]

For $\Psi,\Sigma\in \mathcal{S}_{++}^d$ such that $\Psi\Sigma = \Sigma\Psi$, let the covariance operator $K_{\Psi,\Sigma}:\mathcal{S}^d\to\mathcal{S}^d$ be defined by
\[
K_{\Psi,\Sigma}U \leqdef \frac{\Psi U\Sigma + \Sigma U\Psi}{2}, \qquad U\in\mathcal{S}^d,
\]
and let $\Omega_{\Psi,\Sigma}\in \mathcal{S}_{++}^{q_d}$ be the matrix representation of $K_{\Psi,\Sigma}$ in the Frobenius-orthonormal basis $E_1,\ldots,E_{q_d}$, that is,
\[
(\Omega_{\Psi,\Sigma})_{ab} \leqdef \langle K_{\Psi,\Sigma} E_a, E_b \rangle_F = \tr\{\Sigma E_a\Psi E_b\}, \qquad a,b\in [q_d].
\]
When $\Psi = \Sigma$, this reduces to $K_{\Sigma,\Sigma}U = \Sigma U\Sigma$ and
\[
(\Omega_{\Sigma,\Sigma})_{ab} = \frac{r_a r_b}{2}\{\Sigma_{i_a i_b}\Sigma_{j_a j_b} + \Sigma_{i_a j_b}\Sigma_{j_a i_b}\}, \qquad a,b\in [q_d].
\]

\begin{definition}[Symmetric matrix normal (SN) distribution]\label{def:symmetric.matrix.normal.scale.structured}
Let $M\in \mathcal{S}^d$ and $\Psi,\Sigma\in \mathcal{S}_{++}^d$ be such that $\Psi\Sigma = \Sigma\Psi$. A random matrix $\mathfrak{Z}\in \mathcal{S}^d$ is said to have an SN distribution with mean $M$ and scales $\Psi,\Sigma$, written
\[
\mathfrak{Z}\sim \mathcal{SN}_{d\times d}(M,\Psi,\Sigma),
\]
if $\vecs(\mathfrak{Z})\sim \mathcal{N}_{q_d}(\vecs(M),\Omega_{\Psi,\Sigma})$. In particular, $\EE[\mathfrak{Z}] = M$ and $\Var(\vecs(\mathfrak{Z})) = \Omega_{\Psi,\Sigma}$. The density, with respect to the Lebesgue measure on $\mathcal{S}^d$ induced by the isometry $\vecs$, is
\[
S\mapsto \frac{1}{(2\pi)^{q_d/2}|\Omega_{\Psi,\Sigma}|^{1/2}}\exp\left\{-\frac{1}{2}(\vecs(S-M))^{\top}\Omega_{\Psi,\Sigma}^{-1}\vecs(S-M)\right\}, \qquad S\in \mathcal{S}^d.
\]
\end{definition}

\begin{remark}
If $\Omega_{\Psi,\Sigma}$ is instead parametrized through the unscaled half-vectorization $\vecp(S) \leqdef (S_{i_a j_a})_{a\in [q_d]}\in \R^{q_d}$ used in Section~\ref{sec:Wishart.approximation.shifted.symmetric.normal}, then $\Omega_{\Psi,\Sigma} = R_d \Var(\vecp(\mathfrak{Z})) R_d$, where $R_d$ is the diagonal matrix with diagonal entries $r_a$ such that $\vecs(S) = R_d\vecp(S), ~S\in \mathcal{S}^d$.
\end{remark}

\begin{lemma}[Representation of the SN distribution through symmetrization]\label{lem:symmetrization}
Let $M\in \mathcal{S}^d$ and $\Psi,\Sigma\in \mathcal{S}_{++}^d$ be such that $\Psi\Sigma = \Sigma\Psi$. Let $\mathfrak{G}$ be a $d\times d$ random matrix such that all its entries are independent standard normals. Then
\begin{equation}\label{eq:Z}
\mathfrak{Z} = M + \frac{\Psi^{1/2}\mathfrak{G}\Sigma^{1/2} + \Sigma^{1/2}\mathfrak{G}^{\top}\Psi^{1/2}}{2} \sim \mathcal{SN}_{d\times d}(M,\Psi,\Sigma).
\end{equation}
\end{lemma}

\begin{remark}
The special case $\mathfrak{Z}\sim \mathcal{SN}_{d\times d}(0_{d\times d}, \sqrt{2} I_d, \sqrt{2} I_d)$ is a GOE matrix.
\end{remark}

For any shape parameter $\alpha \in (d-1, \infty)$ and any scale matrix $\Sigma \in \mathcal{S}_{++}^d$, the density of the Wishart distribution, henceforth denoted $\mathcal{W}_d(\alpha, \Sigma)$, is given by
\[
f_{\alpha, \Sigma}^{\mathcal{W}}(X) \leqdef \frac{|X|^{\alpha/2 - (d + 1)/2} \etr(-\Sigma^{-1}X/2)}{|2 \Sigma|^{\alpha/2} \Gamma_d(\alpha/2)}, \qquad X \in \mathcal{S}_{++}^d,
\]
where $\Gamma_d$ denotes the multivariate gamma function.

\newpage
For a time-homogeneous matrix-variate Markov process $(\mathfrak{M}_t)_{t\geq 0}$ taking values in $\mathcal{S}^d$, the transition semigroup of operators $(\mathcal{P}_t)_{t\geq 0}$ is defined, for every measurable function $f$ for which the expectation below is finite, by
\[
\mathcal{P}_t f(Y) = \EE[f(\mathfrak{M}_t) \mid \mathfrak{M}_0 = Y], \qquad t\geq 0.
\]
The corresponding infinitesimal generator of $(\mathfrak{M}_t)_{t\geq 0}$ is defined on its domain by
\[
\mathcal{A} f(Y) = \lim_{s\downarrow 0} \frac{\mathcal{P}_s f(Y) - f(Y)}{s},
\]
provided the limit exists. More generally, for a diffusion, we use the same notation for the extended generator: if $f$ is sufficiently smooth and there exists a measurable function $g$ such that $\smash{(f(\mathfrak{M}_t) - f(\mathfrak{M}_0)-\int_0^t g(\mathfrak{M}_s) \, \rd s)_{t\geq 0}}$ is a local martingale, then we write $\mathcal{A}f = g$. When $f$ belongs to the domain of the infinitesimal generator, the two notions agree.

\section{Stein's method for the SN distribution}\label{sec:symmetric.matrix.normal.Stein.theory}

The symmetric matrix Ornstein--Uhlenbeck (SOU) process considered below is obtained by taking the two-sided matrix Ornstein--Uhlenbeck process of \citet[Eq.~(2)]{GauntOuimetRichards2026} with scale matrices $\Psi$ and $\Sigma$, and replacing the resulting process by its symmetric part around a symmetric mean. Specifically, fix $M\in \mathcal{S}^d$ and $\Psi,\Sigma\in \mathcal{S}_{++}^d$ such that $\Psi\Sigma = \Sigma\Psi$. Consider the $\mathcal{S}^d$-valued process $(\mathfrak{Y}_t)_{t\geq0}$ defined through the following stochastic differential equation (SDE):
\begin{equation}\label{eq:SOU.process}
\rd\mathfrak{Y}_t = - (\mathfrak{Y}_t - M) \, \rd t + \frac{1}{\sqrt{2}} \, \Psi^{1/2}\rd\mathfrak{B}_t\Sigma^{1/2} + \frac{1}{\sqrt{2}} \, \Sigma^{1/2}\rd\mathfrak{B}_t^{\top}\Psi^{1/2}, \qquad \mathfrak{Y}_0 \leqdef Y,
\end{equation}
where $Y\in\mathcal{S}^d$ is deterministic and $(\mathfrak{B}_t)_{t\geq0}$ is a $d\times d$ matrix of independent standard Brownian motions.

For $f\in C^2(\mathcal{S}^d)$ and $S\in\mathcal{S}^d$, define the second-order differential expression
\begin{equation}\label{eq:Frechet.contraction.Psi.nabla.Sigma.nabla}
\begin{aligned}
&\tr\{\Psi\nabla\Sigma\nabla f(S)\} \\
&\leqdef \frac{1}{4}\sum_{i,j=1}^d D^2 f(S)\big[\Psi^{1/2}\bb{e}_i\bb{e}_j^{\top}\Sigma^{1/2} + \Sigma^{1/2}\bb{e}_j\bb{e}_i^{\top}\Psi^{1/2}, \Psi^{1/2}\bb{e}_i\bb{e}_j^{\top}\Sigma^{1/2} + \Sigma^{1/2}\bb{e}_j\bb{e}_i^{\top}\Psi^{1/2}\big].
\end{aligned}
\end{equation}
The analogous expression $\tr\{A\nabla\Sigma\nabla f(S)\}$, for $A\in\mathcal{S}^d$, is understood by linear extension in the first matrix argument. In particular, if $A = \sum_{r=1}^d \eta_r\bb{u}_r\bb{u}_r^{\top}$ is a spectral decomposition, then
\begin{equation}\label{eq:Frechet.contraction.A.nabla.Sigma.nabla}
\tr\{A\nabla\Sigma\nabla f(S)\} = \frac{1}{4}\sum_{r=1}^d \eta_r\sum_{i,j=1}^d \Sigma_{ij} D^2 f(S)[\bb{u}_r\bb{e}_j^{\top} + \bb{e}_j\bb{u}_r^{\top}, \bb{u}_r\bb{e}_i^{\top} + \bb{e}_i\bb{u}_r^{\top}].
\end{equation}

The explicit expression for the extended generator of $(\mathfrak{Y}_t)_{t\geq0}$, denoted $\mathcal{A}_{M,\Psi,\Sigma}^{\mathrm{SOU}}$, is derived in Proposition~\ref{prop:SOU.generator} below.

\begin{proposition}[Extended generator]\label{prop:SOU.generator}
For any $f\in C^2(\mathcal{S}^d)$, we have
\begin{equation}\label{eq:SOU.generator}
\mathcal{A}_{M,\Psi,\Sigma}^{\mathrm{SOU}}f(S) = \tr\{(M-S)\nabla f(S)\} + \tr\{\Psi\nabla\Sigma\nabla f(S)\}, \qquad S\in \mathcal{S}^d.
\end{equation}
\end{proposition}

\begin{proposition}\label{prop:SOU.distribution}
For the SOU process defined in \eqref{eq:SOU.process}, we have
\begin{equation}\label{eq:transition.law}
\mathfrak{Y}_t \mid \{\mathfrak{Y}_0 = Y\} \sim \mathcal{SN}_{d\times d}\big(M + e^{-t}(Y-M),\sqrt{1-e^{-2t}}\,\Psi,\sqrt{1-e^{-2t}}\,\Sigma\big), \qquad t > 0,
\end{equation}
and
\begin{equation}\label{eq:stationary.limiting.distribution}
\mathfrak{Y}_{\infty} \sim \mathcal{SN}_{d\times d}(M,\Psi,\Sigma).
\end{equation}
In particular, let $(\mathcal{P}^{\mathrm{SOU}}_t)_{t\geq0}$ be the transition semigroup with kernel $P_t(Y, \, \rd Z)$, so that, for every bounded Borel measurable function $h:\mathcal{S}^d\to\R$,
\[
(\mathcal{P}^{\mathrm{SOU}}_t h)(Y) \leqdef \int_{\mathcal{S}^d}h(Z)P_t(Y, \, \rd Z) = \EE\big[h(\mathfrak{Y}_t)\mid \mathfrak{Y}_0 = Y\big].
\]
For $\mathfrak{G}$ a $d\times d$ random matrix such that all its entries are independent standard normals, the above shows
\begin{equation}\label{eq:SOU.semigroup.representation}
\mathcal{P}^{\mathrm{SOU}}_t h(Y) = \EE\bigg[h\bigg(M + e^{-t}(Y-M) + \sqrt{1-e^{-2t}} \, \frac{\Psi^{1/2}\mathfrak{G}\Sigma^{1/2} + \Sigma^{1/2}\mathfrak{G}^{\top}\Psi^{1/2}}{2}\bigg)\bigg].
\end{equation}
Moreover, given a probability measure $\mu$ on $\mathcal{S}^d$, the pushed-forward measure $\mu\mathcal{P}^{\mathrm{SOU}}_t$ is defined by
\[
(\mu\mathcal{P}^{\mathrm{SOU}}_t)(A) \leqdef \int_{\mathcal{S}^d}P_t(Y,A) \, \mu(\rd Y), \qquad A\subseteq\mathcal{S}^d \text{ Borel}.
\]
Hence, for $\gamma_{M,\Psi,\Sigma} \leqdef \mathcal{SN}_{d\times d}(M,\Psi,\Sigma)$, we have the invariance
\begin{equation}\label{eq:invariance.SOU}
\gamma_{M,\Psi,\Sigma}\mathcal{P}^{\mathrm{SOU}}_t = \gamma_{M,\Psi,\Sigma}.
\end{equation}
\end{proposition}

This leads to the following Stein characterization for the SN distribution.

\begin{corollary}[Stein characterization]\label{cor:Stein.symmetric.normal}
Let $M\in \mathcal{S}^d$ and $\Psi,\Sigma\in \mathcal{S}_{++}^d$ such that $\Psi\Sigma = \Sigma\Psi$. Then
\[
\mathfrak{X}\sim \mathcal{SN}_{d\times d}(M,\Psi,\Sigma) \qquad \Leftrightarrow \qquad \EE\big[\mathcal{A}_{M,\Psi,\Sigma}^{\mathrm{SOU}}f(\mathfrak{X})\big] = 0 ~~\forall f\in C_{\mathcal{A}_{M,\Psi,\Sigma}^{\mathrm{SOU}}}^2(\mathcal{S}^d),
\]
where, for $\mathfrak{Z}\sim \mathcal{SN}_{d\times d}(M,\Psi,\Sigma)$,
\[
C_{\mathcal{A}_{M,\Psi,\Sigma}^{\mathrm{SOU}}}^2(\mathcal{S}^d)
\leqdef \Big\{f\in C^2(\mathcal{S}^d) \, : \, \EE\big[|\tr\{(M-\mathfrak{Z})\nabla f(\mathfrak{Z})\}|\big] < \infty, ~\EE\big[|\tr\{\Psi\nabla\Sigma\nabla f(\mathfrak{Z})\}|\big] < \infty \Big\}.
\]
\end{corollary}

For any $\beta\in(0,1]$ and $h:\mathcal{S}^d\to\R$, define the $\beta$-H\"older seminorm and the corresponding space of $\beta$-H\"older continuous functions on $\mathcal{S}^d$:
\begin{equation}\label{eq:seminorm.symmetric.normal}
[h]_{\beta,\mathcal{S}} \leqdef \sup_{S\neq T}\frac{|h(S)-h(T)|}{\|S-T\|_F^{\beta}}, \qquad C^{0,\beta}(\mathcal{S}^d) \leqdef \{h:\mathcal{S}^d\to\R \mid [h]_{\beta,\mathcal{S}} < \infty\}.
\end{equation}
The class $C^{0,1}(\mathcal{S}^d)$ is the space of Lipschitz continuous functions on $\mathcal{S}^d$ and $[h]_{1,\mathcal{S}}$ is the minimum Lipschitz constant of $h$. For $p\in\N_0$, we let $\mathrm{Lip}_{p}^{\mathrm{s}}(\mathcal{S}^d)$ denote the class of functions $f$ on $\mathcal{S}^d$ whose coordinate directional derivatives $D_{\bb{a}}^{\mathrm{s}}f$, $\bb{a}\in [q_d]^k$, up to order $k = p$ exist, with the convention that the zeroth-order derivative of a function is the function itself, and whose coordinate directional derivatives of order $p$ are in the class $C^{0,1}(\mathcal{S}^d)$. If $h\in \mathrm{Lip}_{p}^{\mathrm{s}}(\mathcal{S}^d)$ and $\bb{b} = (b_1,\ldots,b_{p + 1})\in [q_d]^{p + 1}$, the coordinate directional derivative $D_{\bb{b}}^{\mathrm{s}} \hspace{0.5mm} h$ is understood to exist almost everywhere, and we write $\|D_{\bb{b}}^{\mathrm{s}} \hspace{0.5mm} h\|_{\infty} \leqdef \operatorname*{ess\,sup}_{S\in\mathcal{S}^d}|D_{\bb{b}}^{\mathrm{s}} \hspace{0.5mm} h(S)|$. For probability measures $\mu$ and $\nu$ on $\mathcal{S}^d$ with finite $\beta$-th moments, let $d_{\mathrm{HK},\beta}^{\mathcal{S}}$ denote the $\beta$-H\"older-Kantorovich distance on $\mathcal{S}^d$ induced by the $\beta$-H\"older seminorm in \eqref{eq:seminorm.symmetric.normal}, that is,
\[
d_{\mathrm{HK},\beta}^{\mathcal{S}}(\mu,\nu) \leqdef \sup\left\{\Big|\int h \, \rd\mu-\int h \, \rd\nu\Big|:h\in C^{0,\beta}(\mathcal{S}^d), ~[h]_{\beta,\mathcal{S}}\leq1\right\}.
\]
Theorem~\ref{thm:Stein.solutions.symmetric.normal} below provides an explicit solution $f_h:\mathcal{S}^d\to\R$ to the SN Stein equation
\begin{equation}\label{eq:Stein.equation.symmetric.normal}
\mathcal{A}_{M,\Psi,\Sigma}^{\mathrm{SOU}}f_h(Y) = h(Y)-\EE[h(\mathfrak{Y}_{\infty})],
\end{equation}
for test functions $h$ belonging either to $C^{0,\beta}(\mathcal{S}^d)$ for some $\beta\in(0,1]$ or to $\mathrm{Lip}_{p}^{\mathrm{s}}(\mathcal{S}^d)$ for some $p\in\N_0$. For $\beta$-H\"older continuous test functions, it also provides a pointwise bound on the solution.

\begin{theorem}[Solution of the SN Stein equation]\label{thm:Stein.solutions.symmetric.normal}
Let $(\mathfrak{Y}_t)_{t\geq0}$ be the SOU process in \eqref{eq:SOU.process} with transition semigroup $(\mathcal{P}^{\mathrm{SOU}}_t)_{t\geq0}$, extended generator $\mathcal{A}_{M,\Psi,\Sigma}^{\mathrm{SOU}}$ given by \eqref{eq:SOU.generator}, and stationary limiting distribution $\gamma_{M,\Psi,\Sigma}$ from Proposition~\ref{prop:SOU.distribution}. For every test function $h$ belonging either to $C^{0,\beta}(\mathcal{S}^d)$ for some $\beta\in(0,1]$ or to $\mathrm{Lip}_{p}^{\mathrm{s}}(\mathcal{S}^d)$ for some $p\in\N_0$, the function
\begin{equation}\label{eq:fh.def.SOU.beta}
f_h(Y) \leqdef -\int_0^{\infty}\Big\{\mathcal{P}^{\mathrm{SOU}}_t h(Y)-\EE[h(\mathfrak{Y}_{\infty})]\Big\} \, \rd t, \qquad Y\in\mathcal{S}^d,
\end{equation}
is well defined pointwise and solves the SN Stein equation \eqref{eq:Stein.equation.symmetric.normal} in the pointwise sense, i.e., for every $Y\in\mathcal{S}^d$, $\lim_{s\downarrow0} \{\mathcal{P}_s^{\mathrm{SOU}}f_h(Y)-f_h(Y)\}/s = h(Y)-\EE[h(\mathfrak{Y}_{\infty})]$. Moreover, if $h\in C^{0,\beta}(\mathcal{S}^d)$ for some $\beta\in(0,1]$, then
\begin{equation}\label{eq:fh.bound.SOU.beta}
|f_h(Y)|
\leq \frac{1}{\beta} [h]_{\beta,\mathcal{S}} \, d_{\mathrm{HK},\beta}^{\mathcal{S}}(\delta_Y,\gamma_{M,\Psi,\Sigma})
\leq \frac{1}{\beta}[h]_{\beta,\mathcal{S}} \Big\{\EE\big[\|\mathfrak{Y}_{\infty} - M\|_F^{\beta}\big] + \|Y - M\|_F^{\beta}\Big\},
\end{equation}
where $\delta_Y$ is the unit point mass at $Y$ and $\mathfrak{Y}_{\infty}\sim \gamma_{M,\Psi,\Sigma}$.
\end{theorem}

Next, we state regularity estimates for the solution to the SN Stein equation.

\begin{theorem}[Regularity of the solution of the SN Stein equation]\label{thm:smoothness.estimates.symmetric.normal}
Let $M\in \mathcal{S}^d$ and $\Psi,\Sigma\in \mathcal{S}_{++}^d$ such that $\Psi\Sigma = \Sigma\Psi$, and let $\mathfrak{Y}_{\infty}\sim \mathcal{SN}_{d\times d}(M,\Psi,\Sigma)$. Let $\bb{a} = (a_1,\ldots,a_m)\in [q_d]^m$. If $h\in \mathrm{Lip}_{m-1}^{\mathrm{s}}(\mathcal{S}^d)$ for some $m\geq1$, then
\begin{equation}\label{eq:Stein.bound.symmetric.normal.1}
\left\|D_{\bb{a}}^{\mathrm{s}} f_h\right\|_{\infty} \leq \frac{1}{m}\left\|D_{\bb{a}}^{\mathrm{s}} h\right\|_{\infty}.
\end{equation}
Alternatively, if $h\in C^{0,\beta}(\mathcal{S}^d)$ is bounded for some $\beta\in(0,1]$, then
\begin{equation}\label{eq:Stein.bound.symmetric.normal.2}
\left\|D_a^{\mathrm{s}} f_h\right\|_{\infty} \leq \sqrt{\dfrac{\pi}{2}}\sqrt{(\Omega_{\Psi,\Sigma}^{-1})_{aa}} ~\|h-\EE[h(\mathfrak{Y}_{\infty})]\|_{\infty}, \qquad a\in [q_d],
\end{equation}
and if $h\in \mathrm{Lip}_{m-2}^{\mathrm{s}}(\mathcal{S}^d)$ with $m\geq2$, then, with $\bb{a}_{-k}$ denoting the vector obtained from $\bb{a}$ by deleting $a_k$,
\begin{equation}\label{eq:Stein.bound.symmetric.normal.3}
\left\|D_{\bb{a}}^{\mathrm{s}} f_h\right\|_{\infty} \leq \frac{\Gamma(m/2)}{\sqrt{2} \, \Gamma(m/2 + 1/2)}\min\limits_{1\leq k\leq m}\left\{\sqrt{(\Omega_{\Psi,\Sigma}^{-1})_{a_k,a_k}} \big\|D_{\bb{a}_{-k}}^{\mathrm{s}}h\big\|_{\infty}\right\}.
\end{equation}
Moreover, for $k\geq1$ and any sufficiently smooth $g:\mathcal{S}^d\to\R$, define the supremum norm in $S$ of the operator norm of the $k$-th Fr\'echet derivative as follows:
\[
\mathcal{M}_k^{\star}(g) \leqdef \sup_{S\in \mathcal{S}^d}\ \sup_{\|U_1\|_F = \cdots = \|U_k\|_F = 1}\big|D^k g(S)[U_1,\ldots,U_k]\big|.
\]
If $h\in C_b^m(\mathcal{S}^d)$ for some $m\geq1$, then
\begin{equation}\label{eq:Stein.bound.symmetric.normal.4}
\mathcal{M}_m^{\star}(f_h) \leq \frac{1}{m}\mathcal{M}_m^{\star}(h).
\end{equation}
If $h\in C^{0,\beta}(\mathcal{S}^d)$ is bounded for some $\beta\in(0,1]$, then
\begin{equation}\label{eq:Stein.bound.symmetric.normal.5}
\mathcal{M}_1^{\star}(f_h) \leq \sqrt{\dfrac{\pi}{2}} \, \|(\Psi\Sigma)^{-1/2}\|_2\|h-\EE[h(\mathfrak{Y}_{\infty})]\|_{\infty},
\end{equation}
and if $h\in C_b^{m-1}(\mathcal{S}^d)$ with $m\geq2$, then
\begin{equation}\label{eq:Stein.bound.symmetric.normal.6}
\mathcal{M}_m^{\star}(f_h) \leq \frac{\Gamma(m/2)}{\sqrt{2} \, \Gamma(m/2 + 1/2)}\|(\Psi\Sigma)^{-1/2}\|_2\mathcal{M}_{m-1}^{\star}(h).
\end{equation}
\end{theorem}

\section{An SN approximation to the Wishart distribution}\label{sec:Wishart.approximation.shifted.symmetric.normal}

When the shape parameter is large, a Wishart matrix behaves like a Gaussian perturbation of its mean. Indeed, when $\alpha\in [d,\infty) \cap \N$, the Gaussian representation of the Wishart distribution \citep[Theorem~10.3.2]{Muirhead1982} shows that $\mathfrak{W}\sim \mathcal{W}_d(\alpha,\Sigma)$ can be written as
\[
\mathfrak{W} \stackrel{\mathrm{law}}{ = } \sum_{r=1}^{\alpha} \bb{X}_r \bb{X}_r^{\top},
\]
where $\bb{X}_1,\ldots,\bb{X}_{\alpha}$ are independent and $\mathcal{N}_d(\bb{0}_d,\Sigma)$-distributed. It is therefore natural to compare $\mathfrak{W}$ with a symmetric Gaussian matrix having the same first two moments. Let $\mathfrak{G}$ be a $d\times d$ matrix of independent standard normal random variables and set $\Sigma_{\alpha} \leqdef \sqrt{2\alpha}\,\Sigma$. By Lemma~\ref{lem:symmetrization},
\begin{equation}\label{eq:Z.alpha.Sigma}
\mathfrak{Z}_{\alpha,\Sigma}
\leqdef \alpha \Sigma + \sqrt{2\alpha} \, \Sigma^{1/2} \frac{\mathfrak{G} + \mathfrak{G}^{\top}}{2} \Sigma^{1/2}
\sim \mathcal{SN}_{d\times d}(\alpha\Sigma,\Sigma_{\alpha},\Sigma_{\alpha}).
\end{equation}
In particular, relative to the Lebesgue measure induced by the unscaled half-vectorization $\vecp$, the density of $\mathfrak{Z}_{\alpha,\Sigma}$ is
\[
f_{\mathfrak{Z}_{\alpha,\Sigma}}(X) = \frac{1}{\sqrt{2^d \pi^{d(d + 1)/2} |\sqrt{2\alpha}\Sigma|^{d + 1}}} \etr\left\{-\frac{1}{4\alpha}\Sigma^{-1}(X - \alpha \Sigma)\Sigma^{-1}(X - \alpha \Sigma)\right\}, \quad X\in \mathcal{S}^d;
\]
see, e.g., \citet[Eq.~(2.5.8)]{GuptaNagar2000}. Moreover, its expectation and covariances are
\[
\EE[\mathfrak{Z}_{\alpha,\Sigma}] = \alpha \Sigma, \qquad \Cov((\mathfrak{Z}_{\alpha,\Sigma})_{ij}, (\mathfrak{Z}_{\alpha,\Sigma})_{k\ell}) = \alpha \{\Sigma_{ik}\Sigma_{j\ell} + \Sigma_{i\ell}\Sigma_{jk}\};
\]
see, e.g., \citet[Theorem~2.5.1]{GuptaNagar2000}. The same expectation and covariance expressions hold for $\mathfrak{W}$; see, e.g., \citet[Theorem~3.3.15]{GuptaNagar2000}.

A local limit theorem for this moment-matched approximation was obtained by \citet{MR4358612}, who derived an asymptotic expansion for the log-ratio of the Wishart density to the corresponding SN density and a resulting total variation bound. The next proposition quantifies this Gaussian approximation by comparing the Wishart extended generator
\begin{equation}\label{eq:Wishart.process.generator}
\mathcal{A}^{\mathcal{W}}f(S) = 2 \, \tr\{(\alpha \Sigma - S) \nabla f(S)\} + 4 \, \tr\{S \nabla \Sigma \nabla f(S)\},
\end{equation}
\citep[see][Proposition~3.1]{BaillyGauntOuimetRichardsvonSachs2026} with the extended generator $\mathcal{A}^{\mathrm{SOU}}$ of an SOU process whose stationary distribution is the law of $\mathfrak{Z}_{\alpha,\Sigma}$.
 Recall from Section~\ref{sec:symmetric.matrix.normal.Stein.theory} that $C^{0,1}(\mathcal{S}^d)$ is the class of Lipschitz continuous functions on $(\mathcal{S}^d,\|\cdot\|_F)$ and that $[h]_{1,\mathcal{S}}$ denotes the minimum Lipschitz constant of $h$. For probability measures $\mu$ and $\nu$ on $\mathcal{S}^d$ with finite first moments, define the Wasserstein distance on $(\mathcal{S}^d,\|\cdot\|_F)$ by
\[
d_{\mathrm{W}}^{\mathcal{S}}(\mu,\nu) \leqdef \sup\left\{\Big|\int h \, \rd\mu-\int h \, \rd\nu\Big|:h\in C^{0,1}(\mathcal{S}^d), \, [h]_{1,\mathcal{S}}\leq1\right\}.
\]
In the notation of Section~\ref{sec:symmetric.matrix.normal.Stein.theory}, $d_{\mathrm{W}}^{\mathcal{S}} = d_{\mathrm{HK},1}^{\mathcal{S}}$.

\begin{proposition}[Wasserstein SN approximation of the Wishart law]\label{prop:Wishart.vs.shifted.symmetric.normal}
Let $\alpha\in (d-1,\infty)$ and $\Sigma\in \mathcal{S}_{++}^d$. Let $\mathfrak{W}\sim \mathcal{W}_d(\alpha,\Sigma)$, and let $\mathfrak{Z}_{\alpha,\Sigma}$ be defined as in \eqref{eq:Z.alpha.Sigma}. Let $\mu_{\mathfrak{W}}$ and $\mu_{\mathfrak{Z}_{\alpha,\Sigma}}$ denote the laws of $(\mathfrak{W} - \alpha \Sigma)/\sqrt{\alpha}$ and $(\mathfrak{Z}_{\alpha,\Sigma} - \alpha \Sigma)/\sqrt{\alpha}$, respectively, on $\mathcal{S}^d$. Then
\begin{equation}\label{eq:Wishart.shifted.symmetric.normal.bound}
d_{\mathrm{W}}^{\mathcal{S}}(\mu_{\mathfrak{W}},\mu_{\mathfrak{Z}_{\alpha,\Sigma}}) \leq \frac{1}{\sqrt{\pi \alpha}} \, \|\Sigma^{-1}\|_2 \, \sqrt{\{(d + 2)\tr(\Sigma)^2 + \|\Sigma\|_F^2\} \{\tr(\Sigma)^2 + \|\Sigma\|_F^2\}}.
\end{equation}
The $\alpha^{-1/2}$ rate in the bound (\ref{eq:Wishart.shifted.symmetric.normal.bound}) is optimal.
\end{proposition}

\begin{remark}
To gauge the dimensional dependence of the bound (\ref{eq:Wishart.shifted.symmetric.normal.bound}), consider the case $\Sigma=I_d$. Then $\|I_d^{-1}\|_2=1$, $\tr(I_d)=d$ and $\|I_d\|_F^2=d$, and hence $d_{\mathrm{W}}^{\mathcal{S}}(\mu_{\mathfrak{W}},\mu_{\mathfrak{Z}_{\alpha,I_d}}) \leq d(d+1)^{3/2}/\sqrt{\pi \alpha}$. In comparison, \citet[Theorem~6]{MR4358612} provides a bound on the rate of convergence in total variation distance: $d_{\mathrm{TV}}^{\mathcal{S}}(\mu_{\mathfrak{W}},\mu_{\mathfrak{Z}_{\alpha,I_d}})\leq C d^{\hspace{0.3mm}3/2}/\sqrt{\alpha}$, as $\alpha\rightarrow\infty$, for some unspecified universal constant $C>0$.
\end{remark}

\section{Proofs}\label{sec:proofs}

\subsection{Proof of Lemma~\ref{lem:symmetrization}}\label{sec:proof.lem:symmetrization}

The map $\mathcal{L}(A) = M + (\Psi^{1/2} A \Sigma^{1/2} + \Sigma^{1/2} A^{\top} \Psi^{1/2})/2$ is affine, so $\mathfrak{Z} = \mathcal{L}(\mathfrak{G})$ is Gaussian in $\mathcal{S}^d$, with $\EE[\mathfrak{Z}] = M$. Moreover, for $U,V\in\mathcal{S}^d$,
\[
\tr\{U(\mathfrak{Z}-M)\} = \tr\{\Sigma^{1/2} U \Psi^{1/2}\mathfrak{G}\}, \qquad \tr\{V(\mathfrak{Z}-M)\} = \tr\{\Sigma^{1/2} V \Psi^{1/2}\mathfrak{G}\}.
\]
Therefore, using the identity $\Cov(\tr\{A\mathfrak{G}\},\tr\{B\mathfrak{G}\}) = \tr(AB^{\top})$, valid for $A,B\in\R^{d\times d}$,
\[
\Cov(\tr\{U\mathfrak{Z}\},\tr\{V\mathfrak{Z}\}) = \tr\{\Sigma^{1/2}U\Psi^{1/2}(\Sigma^{1/2}V\Psi^{1/2})^{\top}\} = \tr\{\Sigma U\Psi V\} = \langle K_{\Psi,\Sigma}U,V\rangle_F.
\]
The last display identifies the covariance operator of $\mathfrak{Z}$ as $K_{\Psi,\Sigma}$, which corresponds to the covariance matrix $\Omega_{\Psi,\Sigma}$ for $\vecs(\mathfrak{Z})$. It follows that $\mathfrak{Z} \sim \mathcal{SN}_{d\times d}(M,\Psi,\Sigma)$.

\subsection{Proof of Proposition~\ref{prop:SOU.generator}}\label{sec:proof.prop:SOU.generator}

Let $(\mathfrak{X}_t)_{t\geq0}$ solve the two-sided matrix Ornstein--Uhlenbeck SDE
\[
\rd\mathfrak{X}_t = -\mathfrak{X}_t \, \rd t + \sqrt{2} \, \Psi^{1/2}\rd\mathfrak{B}_t\Sigma^{1/2}, \qquad \mathfrak{X}_0 = Y-M.
\]
Then $\mathfrak{Y}_t = M + (\mathfrak{X}_t + \mathfrak{X}_t^{\top})/2$ satisfies \eqref{eq:SOU.process}. Applying It\^o's formula to $f(\mathfrak{Y}_t)$ gives the drift contribution $Df(S)[M-S] = \tr\{(M-S)\nabla f(S)\}$. For each Brownian coordinate $\mathfrak{B}_{ij}$, the diffusion direction is $(1/\sqrt{2}) (\Psi^{1/2}\bb{e}_i\bb{e}_j^{\top}\Sigma^{1/2} + \Sigma^{1/2}\bb{e}_j\bb{e}_i^{\top}\Psi^{1/2})$. Hence, by \eqref{eq:Frechet.contraction.Psi.nabla.Sigma.nabla}, the quadratic variation contribution is $\tr\{\Psi\nabla\Sigma\nabla f(S)\}$. This proves \eqref{eq:SOU.generator}.

\subsection{Proof of Proposition~\ref{prop:SOU.distribution}}\label{sec:proof.prop:SOU.distribution}

Solving \eqref{eq:SOU.process} gives
\[
\mathfrak{Y}_t = M + e^{-t}(Y-M) + \frac{1}{\sqrt{2}}\int_0^t e^{-(t-u)}\{\Psi^{1/2}\rd\mathfrak{B}_u\Sigma^{1/2} + \Sigma^{1/2}\rd\mathfrak{B}_u^{\top}\Psi^{1/2}\}.
\]
Since $\int_0^t e^{-(t-u)}\rd\mathfrak{B}_u$ has the same law as $\sqrt{(1-e^{-2t})/2}\,\mathfrak{G}$, Lemma~\ref{lem:symmetrization} gives the transition law \eqref{eq:transition.law} and the semigroup representation \eqref{eq:SOU.semigroup.representation}. Letting $t\to\infty$ yields the stationary limiting distribution \eqref{eq:stationary.limiting.distribution}, and the invariance \eqref{eq:invariance.SOU} follows immediately from the transition law.

\subsection{Proof of Theorem~\ref{thm:Stein.solutions.symmetric.normal}}\label{sec:proof.thm:Stein.solutions.symmetric.normal}

Let $\bb{m} \leqdef \vecs(M)$ and define $\widetilde{h}_{\bb{m}}:\R^{q_d}\to\R$ by $\widetilde{h}_{\bb{m}}(\bb{x}) \leqdef h(\mathrm{mat}_{\mathrm{s}}(\bb{m} + \bb{x}))$. Since $\vecs$ is an isometry, $\widetilde{h}_{\bb{m}}$ belongs to $C^{0,\beta}(\R^{q_d})$ with the same H\"older seminorm as $h$ whenever $h\in C^{0,\beta}(\mathcal{S}^d)$, and belongs to $\mathrm{Lip}_{p}(\R^{q_d})$ whenever $h\in \mathrm{Lip}_{p}^{\mathrm{s}}(\mathcal{S}^d)$. In the $\vecs$ coordinates, the centered process $\vecs(\mathfrak{Y}_t-M)$ is the matrix Ornstein--Uhlenbeck process of \citet[Eq.~(2)]{GauntOuimetRichards2026} with dimensions $q_d\times1$, row-scale $\Omega_{\Psi,\Sigma}$ and column-scale $1$. Let $(\widetilde{\mathcal{P}}_t)_{t\geq0}$ denote its transition semigroup. If $\smash{F_{\widetilde{h}_{\bb{m}}}}$ denotes the corresponding solution of the matrix normal Stein equation from \citet[Theorem~1]{GauntOuimetRichards2026}, then
\[
\begin{aligned}
F_{\widetilde{h}_{\bb{m}}}(\vecs(Y-M))
&= -\int_0^{\infty}\Big\{\widetilde{\mathcal{P}}_t\widetilde{h}_{\bb{m}}(\vecs(Y-M))-\EE[\widetilde{h}_{\bb{m}}(\vecs(\mathfrak{Y}_{\infty}-M))]\Big\} \, \rd t \\
&= -\int_0^{\infty}\Big\{\mathcal{P}^{\mathrm{SOU}}_t h(Y)-\EE[h(\mathfrak{Y}_{\infty})]\Big\} \, \rd t
= f_h(Y),
\end{aligned}
\]
where the last equality follows from \eqref{eq:fh.def.SOU.beta}. The well-definedness, the fact that $f_h$ solves \eqref{eq:Stein.equation.symmetric.normal} pointwise when $\mathcal{A}_{M,\Psi,\Sigma}^{\mathrm{SOU}}$ is interpreted as the infinitesimal generator, and the first bound in \eqref{eq:fh.bound.SOU.beta} follow from \citet[Theorem~1]{GauntOuimetRichards2026} after this translation. The second bound in \eqref{eq:fh.bound.SOU.beta} follows from
\[
d_{\mathrm{HK},\beta}^{\mathcal{S}}(\delta_Y,\gamma_{M,\Psi,\Sigma})\leq\EE\big[\|Y-\mathfrak{Y}_{\infty}\|_F^{\beta}\big]\leq\|Y-M\|_F^{\beta} + \EE\big[\|\mathfrak{Y}_{\infty}-M\|_F^{\beta}\big],
\]
where we used $(u + v)^{\beta}\leq u^{\beta} + v^{\beta}$ for $u,v\geq0$ and $\beta\in(0,1]$.

\subsection{Proof of Theorem~\ref{thm:smoothness.estimates.symmetric.normal}}\label{sec:proof.thm:smoothness.estimates.symmetric.normal}

Let $\bb{m}$, $\widetilde{h}_{\bb{m}}$, $\smash{F_{\widetilde{h}_{\bb{m}}}}$, and $f_h(Y) = \smash{F_{\widetilde{h}_{\bb{m}}}}(\vecs(Y-M))$ be as in the proof of Theorem~\ref{thm:Stein.solutions.symmetric.normal}. Under this identification, the coordinate directional derivatives on $\mathcal{S}^d$ are exactly the corresponding coordinate partial derivatives on $\R^{q_d}$. Therefore, \eqref{eq:Stein.bound.symmetric.normal.1}--\eqref{eq:Stein.bound.symmetric.normal.3} follow from \citet[Theorem~2, Eq.~(10)--(12)]{GauntOuimetRichards2026}.

Since $\vecs$ is an isometry, the operator norms of the Fr\'echet derivatives on $(\mathcal{S}^d,\|\cdot\|_F)$ agree with the corresponding operator norms on $\R^{q_d}$. It remains to identify the covariance-dependent constant in the Fr\'echet bounds. Since $\Psi\Sigma = \Sigma\Psi$, one can choose an orthonormal basis $\bb{u}_1,\ldots,\bb{u}_d$ of common eigenvectors and positive numbers $\psi_1,\ldots,\psi_d$ and $\sigma_1,\ldots,\sigma_d$ such that $\Psi\bb{u}_i = \psi_i\bb{u}_i$ and $\Sigma\bb{u}_i = \sigma_i\bb{u}_i$. The matrices $\bb{u}_i\bb{u}_i^{\top}$, $i\in[d]$, and $(\bb{u}_i\bb{u}_j^{\top} + \bb{u}_j\bb{u}_i^{\top})/\sqrt{2}$, $1\leq i<j\leq d$, form a Frobenius-orthonormal eigenbasis of $K_{\Psi,\Sigma}$, with respective eigenvalues $\psi_i\sigma_i$ for $i\in[d]$ and $(\sigma_i\psi_j + \psi_i\sigma_j)/2$ for $1\leq i<j\leq d$. For $i<j$, the arithmetic-geometric mean inequality gives
\[
\frac{\sigma_i\psi_j + \psi_i\sigma_j}{2} \geq \sqrt{\psi_i\sigma_i\psi_j\sigma_j} \geq \min_{1\leq r\leq d}\psi_r\sigma_r.
\]
Therefore $\lambda_{\min}(K_{\Psi,\Sigma}) = \min_{1\leq r\leq d}\psi_r\sigma_r = \lambda_{\min}(\Psi\Sigma)$. Since $\Omega_{\Psi,\Sigma}$ is a matrix representation of $K_{\Psi,\Sigma}$ in a Frobenius-orthonormal basis, it follows that $\smash{\|\Omega_{\Psi,\Sigma}^{-1/2}\|_2} = \|(\Psi\Sigma)^{-1/2}\|_2$. Thus the Fr\'echet bounds \eqref{eq:Stein.bound.symmetric.normal.4}-\eqref{eq:Stein.bound.symmetric.normal.6} follow from the corresponding bounds in \citet[Theorem~2, Eq.~(13)--(15)]{GauntOuimetRichards2026}.

\subsection{Proof of Corollary~\ref{cor:Stein.symmetric.normal}}\label{sec:proof.cor:Stein.symmetric.normal}

Write $\mathcal{A} \leqdef \smash{\mathcal{A}_{M,\Psi,\Sigma}^{\mathrm{SOU}}}$ for simplicity. We prove the implication from left to right using integration by parts. Let $\mathfrak{X}\sim \mathcal{SN}_{d\times d}(M,\Psi,\Sigma)$, $f\in C_{\mathcal{A}}^2(\mathcal{S}^d)$, $\bb{m} \leqdef \vecs(M)$, $\Omega \leqdef \Omega_{\Psi,\Sigma}$, and let $\varphi$ denote the density of $\mathcal{N}_{q_d}(\bb{m},\Omega)$. Recall that $\smash{\widetilde{f}(\bb{x})} \leqdef f(\mathrm{mat}_{\mathrm{s}}(\bb{x}))$, for $\bb{x}\in \R^{q_d}$. Hence, in the $\vecs$ coordinates,
\[
(\mathcal{A}f)(\mathrm{mat}_{\mathrm{s}}(\bb{x})) = -(\bb{x}-\bb{m})^{\top}\nabla_{\bb{x}}\widetilde{f}(\bb{x}) + \tr\{\Omega\nabla_{\bb{x}}^2\widetilde{f}(\bb{x})\}.
\]
Let $\kappa\in C_c^{\infty}(\R)$ satisfy $0\leq\kappa\leq1$, $\kappa=1$ on $[0,1]$ and $\kappa=0$ on $[2,\infty)$, and set
\[
\kappa_R(\bb{x}) \leqdef \kappa\left(\frac{(\bb{x}-\bb{m})^{\top}\Omega^{-1}(\bb{x}-\bb{m})}{R^2}\right), \qquad R>0.
\]
Since
\[
\begin{aligned}
\big[(\mathcal{A}f)(\mathrm{mat}_{\mathrm{s}}(\bb{x}))\big] \varphi(\bb{x})
&= \big[-(\bb{x}-\bb{m})^{\top}\nabla_{\bb{x}}\widetilde{f}(\bb{x})+\tr\{\Omega\nabla_{\bb{x}}^2\widetilde{f}(\bb{x})\}\big] \varphi(\bb{x}) \\
&= \{\Omega\nabla_{\bb{x}}\widetilde{f}(\bb{x})\}^{\top}\nabla_{\bb{x}}\varphi(\bb{x}) + \tr\{\Omega\nabla_{\bb{x}}^2\widetilde{f}(\bb{x})\}\varphi(\bb{x}) \\
&= \mathrm{div}_{\bb{x}}\{\Omega\nabla_{\bb{x}}\widetilde{f}(\bb{x})\varphi(\bb{x})\},
\end{aligned}
\]
integration by parts gives
\[
\begin{aligned}
&\int_{\R^{q_d}}\kappa_R(\bb{x}) \big[(\mathcal{A}f)(\mathrm{mat}_{\mathrm{s}}(\bb{x}))\big] \varphi(\bb{x}) \, \rd\bb{x} \\
&\quad= -\frac{2}{R^2}\int_{\R^{q_d}}\kappa'\left(\frac{(\bb{x}-\bb{m})^{\top}\Omega^{-1}(\bb{x}-\bb{m})}{R^2}\right)(\bb{x}-\bb{m})^{\top}\nabla_{\bb{x}}\widetilde{f}(\bb{x})\varphi(\bb{x}) \, \rd\bb{x}.
\end{aligned}
\]
The two terms in $\mathcal{A}f(\mathfrak{X})$ are integrable by the definition of $C_{\mathcal{A}}^2(\mathcal{S}^d)$, so dominated convergence applies to the left-hand side as $R\to\infty$. The absolute value of the right-hand side is at most $2\|\kappa'\|_{\infty}R^{-2}\EE[|\tr\{(M-\mathfrak{X})\nabla f(\mathfrak{X})\}|]$, which tends to zero. Therefore $\EE[\mathcal{A}f(\mathfrak{X})] = 0$.

We prove the converse. Assume that the Stein identities in the statement hold for a given $\mathcal{S}^d$-valued random matrix $\mathfrak{X}$. Let $\mathfrak{Y}_{\infty}\sim \mathcal{SN}_{d\times d}(M,\Psi,\Sigma)$, and let $h\in C_c^{\infty}(\mathcal{S}^d)$. Since $C_c^{\infty}(\mathcal{S}^d)\subseteq C^{0,1}(\mathcal{S}^d)\cap C_b^2(\mathcal{S}^d)$, Theorems~\ref{thm:Stein.solutions.symmetric.normal} and \ref{thm:smoothness.estimates.symmetric.normal} apply to $h$ and give the semigroup solution $f_h$ defined in \eqref{eq:fh.def.SOU.beta}, which satisfies $\lim_{s\downarrow0}\{\mathcal{P}^{\mathrm{SOU}}_s f_h(Y)-f_h(Y)\}/s = h(Y)-\EE[h(\mathfrak{Y}_{\infty})]$ for $Y\in \mathcal{S}^d$. Moreover, \eqref{eq:Stein.bound.symmetric.normal.4}, applied with $m = 1$ and $m = 2$, shows that $\mathcal{M}_1^{\star}(f_h) < \infty$ and $\mathcal{M}_2^{\star}(f_h) < \infty$. By Proposition~\ref{prop:SOU.generator} and the agreement of the infinitesimal and extended generators stated at the end of Section~\ref{sec:definitions}, this pointwise semigroup-generator identity is also valid with $\mathcal{A}$ interpreted through the differential expression in \eqref{eq:SOU.generator}, so that $\mathcal{A} f_h(Y) = h(Y)-\EE[h(\mathfrak{Y}_{\infty})]$ for $Y\in \mathcal{S}^d$. Moreover, the second-order term in $\mathcal{A} f_h$ is bounded. Since $h$ is bounded, the Stein equation also shows that the first-order term $\tr\{(M - \, \cdot \, )\nabla f_h(\,\cdot\,)\}$ is bounded. Thus $f_h\in C_{\mathcal{A}}^2(\mathcal{S}^d)$. By the assumed Stein identity, it follows that $0 = \EE\big[\mathcal{A} f_h(\mathfrak{X})\big] = \EE[h(\mathfrak{X})] - \EE[h(\mathfrak{Y}_{\infty})]$. Therefore $\EE[h(\mathfrak{X})] = \EE[h(\mathfrak{Y}_{\infty})]$ for all $h\in C_c^{\infty}(\mathcal{S}^d)$. Since $\mathcal{S}^d$ is a finite-dimensional Euclidean space under the Frobenius inner product, the class $C_c^{\infty}(\mathcal{S}^d)$ is measure determining for probability measures on $\mathcal{S}^d$. It follows that $\smash{\mathfrak{X} \stackrel{\mathrm{law}}{=} \mathfrak{Y}_{\infty} \sim \mathcal{SN}_{d\times d}(M,\Psi,\Sigma)}$. This concludes the proof.

\subsection{Proof of Proposition~\ref{prop:Wishart.vs.shifted.symmetric.normal}}\label{sec:proof.prop:Wishart.vs.shifted.symmetric.normal}

By Propositions~\ref{prop:SOU.generator} and \ref{prop:SOU.distribution}, applied with mean $M = \alpha\Sigma$ and both scale matrices equal to $\Sigma_{\alpha} \equiv \sqrt{2\alpha}\,\Sigma$, the law of $\mathfrak{Z}_{\alpha,\Sigma}$ from \eqref{eq:Z.alpha.Sigma} is the stationary limiting distribution of the corresponding SOU process, whose extended generator is
\begin{equation}\label{eq:sym.OU.generator.application4.shifted}
\mathcal{A}^{\mathrm{SOU}} f(X) = \tr\{(\alpha \Sigma - X) \nabla f(X)\} + 2\alpha \, \tr\{\Sigma \nabla \Sigma \nabla f(X)\}, \qquad X\in \mathcal{S}^d.
\end{equation}
Fix $h\in C^{0,1}(\mathcal{S}^d)$ with $[h]_{1,\mathcal{S}}\leq 1$. Since subtracting a constant from $h$ does not change the difference $\EE[h(\mathfrak{W})]-\EE[h(\mathfrak{Z}_{\alpha,\Sigma})]$, we may replace $h$ by $h-h(0_{d\times d})$ and hence assume that $h(0_{d\times d})=0$. It then follows that $|h(X)|\leq\|X\|_F$. For $R>0$, let $\pi_R(X) = [\ind_{\{\|X\|_F\leq R\}} + (R / \|X\|_F) \ind_{\{\|X\|_F>R\}}] X$ be the Frobenius-metric projection onto $\{X\in\mathcal{S}^d:\|X\|_F\leq R\}$ and set $h_R \leqdef h\circ\pi_R$; then $h_R$ is bounded, $[h_R]_{1,\mathcal{S}}\leq1$, $h_R\to h$ pointwise, and dominated convergence applies because $|h_R(X)|\leq\|X\|_F$ and $\mathfrak{W}$ and $\mathfrak{Z}_{\alpha,\Sigma}$ have finite first moments. Finally, mollifying $h_R$ in Frobenius-orthonormal coordinates gives functions $h_{R,\varepsilon}\in C_b^{\infty}(\mathcal{S}^d)$ such that $\mathcal{M}_1^{\star}(h_{R,\varepsilon})\leq1$ and $h_{R,\varepsilon}\to h_R$ uniformly. Hence it suffices to prove the estimate for $h\in C_b^1(\mathcal{S}^d)$ with $\mathcal{M}_1^{\star}(h)\leq1$. We assume this from now on. Let $f_h$ be the corresponding solution of the SN Stein equation from Theorem~\ref{thm:Stein.solutions.symmetric.normal}, with $\beta=1$. Then $\mathcal{A}^{\mathrm{SOU}} f_h(Y) = h(Y)-\EE[h(\mathfrak{Z}_{\alpha,\Sigma})]$, for $Y\in \mathcal{S}^d$. Moreover, by the Fr\'echet derivative bound \eqref{eq:Stein.bound.symmetric.normal.6} in Theorem~\ref{thm:smoothness.estimates.symmetric.normal}, applied with both scale matrices equal to $\Sigma_{\alpha}$ and $m = 2$, we have
\[
\mathcal{M}_2^{\star}(f_h)
\leq \frac{\Gamma(1)}{\sqrt{2} \, \Gamma(3/2)}\|\Sigma_{\alpha}^{-1}\|_2\mathcal{M}_1^{\star}(h)
\leq \frac{1}{\sqrt{\pi\alpha}}\|\Sigma^{-1}\|_2.
\]
Therefore, by the definition of $\mathcal{M}_2^{\star}$, for all $X,U,V\in \mathcal{S}^d$,
\begin{equation}\label{eq:sym.OU.second.derivative.bound.application4.shifted}
\big|D^2 f_h(X)[U,V]\big| \leq \frac{1}{\sqrt{\pi\alpha}}\|\Sigma^{-1}\|_2\|U\|_F\|V\|_F.
\end{equation}

We now compare $\mathfrak{W}$ with $\mathfrak{Z}_{\alpha,\Sigma}$. Since $\mathfrak{W}$ has a finite first moment and $f_h$ has bounded first and second Fr\'echet derivatives, the restriction of $f_h$ to $\mathcal{S}_{++}^d$ belongs to the admissible test-function class in the Wishart Stein characterization of \citet[Corollary~3.3]{BaillyGauntOuimetRichardsvonSachs2026}. Applying the forward implication of that result to $\mathfrak{W}\sim \mathcal{W}_d(\alpha,\Sigma)$, with the Wishart extended generator \eqref{eq:Wishart.process.generator}, gives
$\EE[\mathcal{A}^{\mathcal{W}} f_h(\mathfrak{W})] = 0$.
Combining it with the Stein equation for $\mathcal{A}^{\mathrm{SOU}}$ and the extended-generator identities \eqref{eq:Wishart.process.generator} and \eqref{eq:sym.OU.generator.application4.shifted}, we obtain
\begin{equation}\label{eq:Stein.identity.shifted.symmetric.normal}
\begin{aligned}
\big|\EE[h(\mathfrak{W})] - \EE[h(\mathfrak{Z}_{\alpha,\Sigma})]\big|
&= \big|\EE[\mathcal{A}^{\mathrm{SOU}} f_h(\mathfrak{W})] - (1/2) \, \EE[\mathcal{A}^{\mathcal{W}} f_h(\mathfrak{W})]\big| \\
&= 2 \, \big|\EE[\tr\{(\mathfrak{W} - \alpha \Sigma) \nabla \Sigma \nabla f_h(\mathfrak{W})\}]\big| \\
&\leq 2 \, \EE\big[|\tr\{(\mathfrak{W} - \alpha \Sigma) \nabla \Sigma \nabla f_h(\mathfrak{W})\}|\big].
\end{aligned}
\end{equation}

It remains to bound the right-hand side of \eqref{eq:Stein.identity.shifted.symmetric.normal}. Fix $X\in \mathcal{S}_{++}^d$ and write the spectral decompositions
\vspace{-2.8mm}
\[
X - \alpha \Sigma = \sum_{r=1}^d \eta_r \bb{u}_r \bb{u}_r^{\top}, \qquad \Sigma = \sum_{a=1}^d \lambda_a \bb{v}_a \bb{v}_a^{\top},
\]

\newpage
\noindent
where $\eta_1,\ldots,\eta_d\in \R$, $\lambda_1,\ldots,\lambda_d>0$, and $\bb{u}_1,\ldots,\bb{u}_d$ and $\bb{v}_1,\ldots,\bb{v}_d$ are orthonormal bases of eigenvectors. For $r,a\in [d]$, define $H_r^{(a)} \leqdef \bb{u}_r \bb{v}_a^{\top} + \bb{v}_a \bb{u}_r^{\top}\in \mathcal{S}^d$. Then, by the definition \eqref{eq:Frechet.contraction.A.nabla.Sigma.nabla}, the spectral decomposition of $\Sigma$, and the bilinearity of $D^2f_h(X)$,
\[
\begin{aligned}
\bb{u}_r^{\top}\nabla\Sigma\nabla f_h(X)\bb{u}_r
&= \frac{1}{4}\sum_{i,j=1}^d \Sigma_{ij}D^2 f_h(X)\big[\bb{u}_r\bb{e}_j^{\top}+\bb{e}_j\bb{u}_r^{\top},\bb{u}_r\bb{e}_i^{\top}+\bb{e}_i\bb{u}_r^{\top}\big] \\
&= \frac{1}{4}\sum_{a=1}^d\lambda_aD^2 f_h(X)\left[\sum_{j=1}^d(\bb{v}_a)_j(\bb{u}_r\bb{e}_j^{\top}+\bb{e}_j\bb{u}_r^{\top}),\sum_{i=1}^d(\bb{v}_a)_i(\bb{u}_r\bb{e}_i^{\top}+\bb{e}_i\bb{u}_r^{\top})\right] \\
&= \frac{1}{4}\sum_{a=1}^d \lambda_a D^2 f_h(X)[H_r^{(a)},H_r^{(a)}].
\end{aligned}
\]
By \eqref{eq:sym.OU.second.derivative.bound.application4.shifted} and $\|H_r^{(a)}\|_F^2 = 2 \, \|\bb{u}_r\|_2^2\|\bb{v}_a\|_2^2 + 2 \, \tr\{(\bb{u}_r \bb{v}_a^{\top})^{\top} \bb{v}_a \bb{u}_r^{\top}\} = 2 + 2 \, \bb{u}_r^{\top} \bb{v}_a \bb{v}_a^{\top} \bb{u}_r$, we obtain
\[
\big|\bb{u}_r^{\top} \nabla \Sigma \nabla f_h(X) \bb{u}_r\big|
\leq \frac{1}{4\sqrt{\pi\alpha}}\|\Sigma^{-1}\|_2 \sum_{a=1}^d \lambda_a\|H_r^{(a)}\|_F^2
= \frac{1}{2\sqrt{\pi\alpha}}\|\Sigma^{-1}\|_2\{\tr(\Sigma) + \bb{u}_r^{\top}\Sigma\bb{u}_r\}.
\]
Also, using $\sum_{r=1}^d (\bb{u}_r^{\top}\Sigma\bb{u}_r)^2 \leq \sum_{r,s=1}^d (\bb{u}_r^{\top}\Sigma\bb{u}_s)^2 = \|\Sigma\|_F^2$, we have
\[
\sum_{r=1}^d\{\tr(\Sigma) + \bb{u}_r^{\top}\Sigma\bb{u}_r\}^2
\leq (d + 2) \tr(\Sigma)^2 + \|\Sigma\|_F^2
\equiv B_{\Sigma}.
\]
Hence, by the Cauchy--Schwarz inequality and $(\sum_{r=1}^d \eta_r^2)^{1/2} = \|X - \alpha \Sigma\|_F$,
\[
\big|\tr\{(X - \alpha \Sigma) \nabla \Sigma \nabla f_h(X)\}\big|
= \left|\sum_{r=1}^d \eta_r \, \bb{u}_r^{\top} \nabla \Sigma \nabla f_h(X) \bb{u}_r\right|
\leq \frac{1}{2\sqrt{\pi\alpha}} \, \|\Sigma^{-1}\|_2 \, \sqrt{B_{\Sigma}} \, \|X - \alpha \Sigma\|_F.
\]
Applying this bound with $X = \mathfrak{W}$ in \eqref{eq:Stein.identity.shifted.symmetric.normal} gives
\[
\big|\EE[h(\mathfrak{W})] - \EE[h(\mathfrak{Z}_{\alpha,\Sigma})]\big| \leq \frac{1}{\sqrt{\pi\alpha}} \, \|\Sigma^{-1}\|_2 \, \sqrt{B_{\Sigma}} \, \EE[\|\mathfrak{W} - \alpha \Sigma\|_F].
\]
By Cauchy--Schwarz and the covariance formula for the Wishart distribution \citep[Theorem~3.3.15]{GuptaNagar2000},
\vspace{-3mm}
\[
\EE[\|\mathfrak{W} - \alpha \Sigma\|_F]
\leq \sqrt{\EE[\|\mathfrak{W} - \alpha \Sigma\|_F^2]}
= \sqrt{\alpha \sum_{i,j=1}^d (\Sigma_{ii}\Sigma_{jj} + \Sigma_{ij}^2)}
= \sqrt{\alpha} \, \sqrt{\tr(\Sigma)^2 + \|\Sigma\|_F^2}.
\]
Combining the last two displays gives, for every $h\in C_b^1(\mathcal{S}^d)$ with $\mathcal{M}_1^{\star}(h)\leq 1$,
\[
\big|\EE[h(\mathfrak{W})] - \EE[h(\mathfrak{Z}_{\alpha,\Sigma})]\big| \leq \frac{1}{\sqrt{\pi}} \, \|\Sigma^{-1}\|_2 \, \sqrt{\{(d + 2)\tr(\Sigma)^2 + \|\Sigma\|_F^2\} \{\tr(\Sigma)^2 + \|\Sigma\|_F^2\}}.
\]
By the approximation reduction made above, the same bound holds for every $h\in C^{0,1}(\mathcal{S}^d)$ with $[h]_{1,\mathcal{S}}\leq 1$. Taking the supremum over this class of test functions gives the corresponding raw-scale bound; by the translation invariance and homogeneity of $d_{\mathrm{W}}^{\mathcal{S}}$, centering at $\alpha \Sigma$ and dividing by $\sqrt{\alpha}$ yields \eqref{eq:Wishart.shifted.symmetric.normal.bound}.

We now prove that the $\alpha^{-1/2}$ rate is optimal. Let $\alpha>2$. It suffices to consider the univariate $d = 1$ case with $\Sigma = 1$. In this case, let $T_\alpha \leqdef (\mathfrak{W} - \alpha)/\sqrt{2\alpha}$ denote a normalized $\mathrm{Gamma}(\alpha/2, 1/2)$ random variable in the shape-rate parametrization, which is supported on $(-\sqrt{\alpha/2},\infty)$, and let $Z\sim \mathcal{N}(0,1)$. Write $\mu_{T_\alpha}$ and $\mu_Z$ for their laws. A routine calculation shows that, uniformly for $x\in[-1,2]$, the density $f_{T_\alpha}$ of $T_\alpha$ enjoys the Edgeworth expansion $f_{T_\alpha}(x) = (1/\sqrt{2\pi})e^{-x^2/2} [1+ \{\sqrt{2}/(3\sqrt{\alpha})\}(x^3-3x)+\OO(\alpha^{-1})]$, as $\alpha\rightarrow\infty$. Now let $h_{\star}:\R\to[0,1]$ be defined by $h_{\star}(x) \leqdef 1+x$ for $x\in[-1,0)$, $1-x/2$ for $x\in[0,2]$, and $0$ otherwise. Identifying $\R$ with $\mathcal{S}^1$, observe that $h_{\star}\in C^{0,1}(\mathcal{S}^1)$ with $[h_{\star}]_{1,\mathcal{S}}\leq 1$. Given that $(\mathfrak{Z}_{\alpha,1}-\alpha)/\sqrt{\alpha} ~~ \smash{\stackrel{\mathrm{law}}{=}} \sqrt{2} Z$ and $(\mathfrak{W}-\alpha)/\sqrt{\alpha} = \sqrt{2} \, T_\alpha$, the homogeneity

\newpage
\noindent
of the Wasserstein distance gives $d_{\mathrm{W}}^{\mathcal{S}}(\mu_{\mathfrak{W}},\mu_{\mathfrak{Z}_{\alpha,1}}) = \sqrt{2} \, d_{\mathrm{W}}^{\mathcal{S}}(\mu_{T_\alpha},\mu_Z)$, and thus
\[
\begin{aligned}
&d_{\mathrm{W}}^{\mathcal{S}}(\mu_{\mathfrak{W}},\mu_{\mathfrak{Z}_{\alpha,1}})
\geq \sqrt{2} \, |\EE[h_{\star}(T_\alpha)]-\EE[h_{\star}(Z)]| \\
&\quad= \frac{\sqrt{2}}{3\sqrt{\pi\alpha}}\left|\int_{-1}^0(1 + x)(x^3-3x)e^{-x^2/2} \, \rd x + \int_{0}^2(1-x/2)(x^3-3x)e^{-x^2/2} \, \rd x\right| + \OO(\alpha^{-1})\\
&\quad= \frac{\sqrt{2}(1-e^{-3/2})}{3\sqrt{e\pi\alpha}} + \OO(\alpha^{-1}).
\end{aligned}
\]
This completes the proof.

\addcontentsline{toc}{section}{References}

\bibliographystyle{plainnat}
\bibliography{bib_clean}

%\begin{acks}
%***
%\end{acks}

\end{document}